\documentclass[10pt]{article}
\usepackage{amsmath,amsthm,amsfonts,amssymb}
\usepackage{graphicx}
\usepackage{hyperref}
\usepackage[usenames]{color}
\usepackage{epsfig}

\newtheorem{theorem}{Theorem}[section]
\newtheorem{lemma}[theorem]{Lemma}
\theoremstyle{definition}
\newtheorem{algorithm}[theorem]{Algorithm}
\newtheorem{corollary}[theorem]{Corollary}

\newtheorem{example}[theorem]{Example}
\newtheorem{problem}[theorem]{Problem}
\newtheorem{remark}[theorem]{Remark}
\newtheorem{proposition}[theorem]{Proposition}

\def\CC{{\mathcal C}}

\def\CA{\mathcal A}
\def\CM{\mathcal M}
\def\CV{\mathcal V}

\newcounter{comcount}

\def\MZ{{\mathbb{Z}}}
\def\MN{{\mathbb{N}}}
\def\MR{{\mathbb{R}}}

\title{The Word and Geodesic Problems in Free Solvable Groups}

\author{A. Myasnikov, V. Roman'kov, A. Ushakov, A.Vershik}

\begin{document}
\maketitle

\begin{abstract}
We study the computational complexity of the Word Problem (WP) in free
solvable groups $S_{r,d}$, where
$r \geq 2$ is the  rank and $d \geq 2$ is the solvability  class of the group.
It is known that the Magnus embedding of $S_{r,d}$ into matrices
provides a polynomial time decision algorithm for  WP in a fixed group $S_{r,d}$.
Unfortunately, the degree of the  polynomial grows together with  $d$, so the uniform algorithm is not polynomial in $d$. In this paper we show that WP  has time complexity
$O(r n \log_2 n)$ in  $S_{r,2}$,  and $O(n^3 r d)$ in  $S_{r,d}$  for $d  \geq 3$.  However, it turns out,  that a seemingly close problem of computing the geodesic length of elements in $S_{r,2}$ is $NP$-complete. We prove also  that one can compute Fox derivatives of elements from $S_{r,d}$ in time  $O(n^3 r d)$, in particular one can  use efficiently  the Magnus embedding in computations with free solvable groups.  Our approach is based on such classical tools as the Magnus embedding and Fox calculus, as well as, on  a relatively new  geometric ideas, in particular, we  establish a direct link between Fox derivatives and geometric flows on Cayley graphs.
\end{abstract}

\tableofcontents

\section{Introduction}
\label{se:intro}

In this paper we study the computational complexity of several algorithmic problems related to the Word Problem (WP) in  free solvable groups. Let $S_{r,d}$ be a free solvable group of  rank $r \geq 2$ and the solvability class  $d \geq 2$. We present here a uniform decision algorithm that solves WP in  time $O(r n \log_2 n)$ in  the free metabelian group $S_{r,2}$ (also denoted by $M_r$),  and $O(n^3 r d)$ in the free solvable group  $S_{r,d}$  for $d  \geq 3$, where $n$ is the length of the input word. In particular, this algorithm is at most cubic in $n$ and linear in $r$ and $d$ for all free solvable groups $S_{r,d}$.  Notice, that in all previously known polynomial time decision algorithms for WP in $S_{r,d}$ the degree of the polynomial grows together with $d$.  In fact, we prove more, we show that one can compute Fox derivatives of elements from $S_{r,d}$ in time  $O(n^3 r d)$. This allows one to use efficiently  the Magnus embedding in computations with free solvable groups. On the other hand, we describe geodesics in $S_{r,d}$ and show that a seemingly close problem of finding the geodesic length of a given element from $S_{r,2}$ is surprisingly hard -- it is NP-complete. Our approach is based on such classical tools as the Magnus embedding and Fox calculus, as well as, on  a relatively new (in group theory) geometric ideas from \cite{DLS} and \cite{Vershik_Dobrynin:2004}. In particular, we  establish a direct link between Fox derivatives and geometric flows on Cayley graphs.

The study of algorithmic problems in free solvable groups can be traced to the the work \cite{Magnus:1939} of Magnus, who in 1939 introduced an embedding (now called the Magnus embedding)  of an arbitrary  group of the type $F/N^\prime$ into a matrix group of a particular type with coefficients in the group ring of $F/N$ (see section \ref{subsec:Magnus-embedding} below).  Since WP in free abelian groups is decidable in polynomial time, by induction,  this embedding immediately gives a polynomial time decision algorithm for a fixed free solvable group $S_{r,d}$. However the degree of the polynomial here grows together with $d$.

In 1950's R. Fox introduced his free differential calculus and made the Magnus embedding much more transparent \cite{Fox_calc1,Fox_calc2,Fox_calc3,Fox_calc4} (see also Section \ref{subsec:Fox-der}). Namely, besides other things, he showed that an element $w \in F$  belongs to $N^\prime = [N,N]$ if and only if all partial derivatives of $w$ is equal to $0$ in the integer group ring of $F/N$. This reduces WP in $F/N^\prime$ directly to the word problem in $F/N$. In particular, it solves, by induction, WP in $S_{r,d}$. Again, the decision algorithm is polynomial in a fixed  group $S_{r,d}$, but the degree of the polynomial grows with $d$ -- not a surprise since the partial derivatives of $w$ describe precisely the image of $w$ under the Magnus embedding.

A few years later  P. Hall  proved the finiteness
of all subdirect indecomposable finitely generated
abelian-by-nilpotent groups. This implies that all finitely generated
abelian-by-nilpotent, in particular, metabelian,  groups are
residually finite. About the same time Gruenberg  extended this
result to  arbitrary free solvable groups \cite{Gruenberg:1957}.
Now one can solve WP in $S_{r,d}$ in the following way.
Given $w \in S_{r,d}$, as a word in the fixed set of
generators, one can start two processes in parallel.
The first one enumerates effectively all consequences
of the defining relations of $S_{r,d}$ in $F_r$ (which is possible
since the group is recursively presented) until the word
$w$ occurs,  and the second one enumerates all homomorphisms from $S_{r,d}$ into all finite symmetric groups $S_n$ (checking if a given $r$-tuple of elements in $S_n$ generates a solvable group of class $d$) until it finds one where the  image of $w$ is non-trivial.  However, computer experiments show that the algorithm described above is extremely inefficient (though its complexity is unknown).

Another shot at WP in metabelian groups comes from their linear representations. V. Remeslennikov proved in \cite{Remeslennikov:1969} that a finitely
generated metabelian group  (under some restrictions) is
embeddable  into
$GL(n,R)$ for a suitable $n$ and a suitable  ring $R = K_1 \times
\ldots \times K_n$ which is a finite direct product of fields
$K_i$.   In \cite{Wehrfritz:1973}, see also \cite{Wehrfritz:1976},  B. Wehrfritz generalized this result to arbitrary finitely generated metabelian groups $G$.   It follows that $G$ is embeddable into a finite direct product of linear groups. Since WP in linear groups is polynomial time decidable  this implies that WP in $G$ is polynomial time decidable. Notice, that it is unclear if there is a uniform polynomial time decision algorithm for WP in arbitrary finitely generated metabelian groups.

In comparison, observe, that there are finitely presented solvable groups of class 3 with undecidable WP. In \cite{Kharlampovich:1981} O. Kharlampovich constructed the first example of  such a group by simulating a universal Minski machine in WP of the group.  There are several results which clarify the boundary between decidability and undecidability of the  word problems in solvable groups, we refer to a survey \cite{Olga} for details.

Our approach to WP in free solvable groups is based on the Fox Theorem mentioned above. Using binary tree search techniques and associative arrays we were able to compute Fox's derivatives of elements $w$ of a free solvable group $S_{d,r}$ in time  $O(n^3  d)$, where $n = |w|$.  Significance of this result goes beyond WP for these groups  - it gives  a fast algorithm to compute images of elements under the Magnus embedding. This  opens  up an opportunity to solve effectively other algorithmic problems in groups $S_{r,d}$ using the classical techniques developed for wreath products of groups.

In the second half of the paper, Section \ref{se:metabelian_geodesic},   we study algorithmic problems on geodesics in free metabelian groups.
Let $G$ be a group with a finite   set of generators $X = \{x_1, \ldots,
x_r\}$ and $\mu:F(X) \rightarrow G$ the canonical epimorphism.
For a word $w$ in the alphabet $X^{\pm 1}$ by $|w|$ we denote the length of $w$.  The {\em geodesic length} $l_X(g)$ of
an element $g \in G$ relative to $X$ is defined by
$$
l_X(g) = \min \{|w| \mid w \in F(X), w^\mu = g\}.
$$
We write, sometimes,  $l_X(w)$ instead of $l_X(w^\mu)$.  A word $w \in F(X)$ is called {\em geodesic} in $G$ relative to $X$, if $|w| = l_X(w)$.
We are interested here in the following two algorithmic {\em search} problems in  $G$.

\medskip\noindent
{\bf The Geodesic Problem (GP):} \ Given a word $w \in F(X)$ find a
word $u \in F(X)$ which is geodesic in $G$ such that $w^\mu = u^\mu$.

\medskip
\noindent
 {\bf The Geodesic Length Problem (GLP):}
\  Given a word $w \in F(X)$ find $l_X(w)$.

\medskip
Though GLP seems easier than GP, in practice, to solve GLP one usually solves GP first, and only then computes the geodesic length.
It is an interesting question if there exists a group $G$ and a finite set $X$ of generators
for $G$ relative to which GP is strictly harder than GLP.

As customary in complexity theory  one can modify the search problem GLP to get the corresponding bounded  {\em decision} problem (that requires only answers "yes" or "no"):

 \medskip
\noindent
 {\bf   The Bounded Geodesic Length Problem (BGLP):}
\ Let $G$ be a group with a finite generating set $X$. Given a word $w \in F(X)$ and  a natural number $k$ determine if $l_X(w) \leq k$.

\medskip
In Section \ref{subsec:geodesics-in-groups} we compare in detail   the algorithmic "hardness" of the problems WP, BGLP, GLP, and GP in a given group $G$. Here we would like only to mention that  in the list of the problems above each one is Turing reducible in polynomial time to the next one in the list, and GP is Turing reducible to WP in exponential time (see definitions in Section \ref{subsec:geodesics-in-groups}).

Among general facts on   computational complexity of geodesics notice that  if $G$ has  a polynomial {\em growth}, i.e., there is a polynomial $p(n)$ such that for each $n \in  \MN$ cardinality of the ball $B_n$  of radius $n$ in the Cayley graph $\Gamma(G,X)$ is at most $p(n)$, then one can easily construct this ball $B_n$ in polynomial time with an oracle for WP in $G$. If, in addition,  such a group  $G$ has WP decidable in polynomial time then all the problems above  have polynomial time complexity with respect to any finite generating set of $G$ (since the growth and WP stay polynomial for any finite set of generators). Now, by Gromov's theorem \cite{Gromov_pgrowth:1981},  groups of polynomial growth are virtually nilpotent, hence linear, so they have WP decidable in polynomial time.  It follows that  all  Geodesic Problems are polynomial time decidable in groups of polynomial growth (finitely generated virtually nilpotent groups).  On the other hand, there are many groups of exponential growth where GP is  decidable in polynomial time, for example, hyperbolic  groups \cite{Epstein} or metabelian Baumslag-Solitar group $BS(1,n) =  \langle a,t \mid t^{-1} a t = a^n\rangle$, $n \geq 2$ (see \cite{Elder_BS:2007} and Section \ref{subsec:geodesics-in-groups} for comments).

In general, if WP in $G$ is  decidable in polynomial time then $BGLP$ is in the class  NP, i.e., it is decidable in polynomial time by a non-deterministic Turing machine.
It might happen though,  that BGLP in a group $G$ is as hard as any in the class NP, i.e., it is NP-complete.  The simplest  example of this type is due to   Perry, who  showed in \cite{Parry:1992} that BGLP  is NP-complete in the metabelian group $\MZ_2 wr (\MZ \times \MZ)$ (the wreath product of $\MZ_2$ and $\MZ \times \MZ$).  Correspondingly, the search  problems GP and GLP are {\em NP-hard}, which means, precisely, that some (any) NP-complete problem is Turing reducible to them in polynomial time.

Our view-point on geodesics  in free solvable groups is based on geometric ideas from the following two papers. In 1993 Droms, Lewin, and Servatius introduced a new  geometric approach to study WP and GLP in  groups of the type $F/N^\prime$   via paths in the Cayley graph of $F/N$ \cite{DLS}.
In 2004 Vershik and Dobrynin studied algebraic structure of solvable  groups, using homology of related Cayley graphs  \cite{Vershik_Dobrynin:2004}. This approach was outlined earlier in the papers \cite{Vershik:1999,Vershik:2000}, where possible applications to random walks on metabelian groups have been discussed. In the papers \cite{Vershik:1999,Vershik:2000}  (see also \cite{Vershik_Dobrynin:2004}) a new robust presentation of a free metabelian group $S_{r,2}$  was introduced as  an extension of $\mathbb{Z}^r$ by the integer first homology group of the lattice $\mathbb{Z}^r$ (viewed as a one-dimensional complex) with a distinguished  2-cocycle. Similar presentations of other metabelian and solvable groups laid out foundations of a new approach to algorithmic problems in solvable groups.

It seems these ideas are still underdeveloped in group-theoretic context, despite their obvious potential. Meanwhile, in semigroup theory similar geometric techniques have been widely used to deal with   free objects in semidirect products of varieties. One can find an explicit exposition of these techniques in the papers due to Almeida \cite{Almeida:1989} and Almeida and Weil \cite{Almeida-Weil:1995}, while in \cite{Rhodes-Steinberg:2001,Auinger-Steinberg:2005,Auinger-Steinberg:2005B,Auinger-Steinberg:2004}
Auinger, Rhodes and Steinberg use similar machinery  on a regular basis.  Earlier, similar methods, though sometimes implicitly, were used in inverse semigroups theory, we refer here to papers \cite{Munn:74,Margolis-Meakin:1989,Margolis-Meakin-Stephen:1990,Cowan:1991}.

In group theory most of the results in this area  relied on various forms of the Magnus embedding  and Fox derivatives. The role that  the Magnus embeddings play in varieties of groups was clarified by  Shmelkin \cite{Shmelkin_Magnus_embed}.   In \cite{Matthews:1966} Matthews  proved that the conjugacy problem (CP) in free metabelian groups is decidable, and Kargapolov and Remeslennikov generalized this to free solvable groups $S_{r,d}$ \cite{Kargapolov-Rem:1966}. A few years later
Remeslennikov and Sokolov described precisely the image of $F/N^\prime$ under the Magnus embedding and showed  that CP is residually finite in  $S_{r,d}$ \cite{Remeslennikov-Sokolov:1970}.   We refer to a survey \cite{Rem-Rom:1980} on algorithmic problems in solvable groups.

In Sections \ref{subsec:flows} and \ref{subsec:geom-interpret} we study elements of groups of the type $F/N^\prime$ via flows on the Cayley graph $\Gamma$ of $F/N$. It turns out the flow generated by a word $w \in F$ on the graph $\Gamma$ directly corresponds to the Fox derivatives of $w$ in the group ring $\MZ F/N$. This simple observation links together the techniques developed in group theory for the Magnus embeddings  with the extensive geometric and the graph-theoretic machinery for flows. Indeed, the set of geometric circulations (flows where the  Kirchhoff law holds for all vertices, including the source and the sink)  form a group which is naturally isomorphic to the first homology group $H_1(\Gamma, \MZ)$ of $\Gamma$. In this content the geometric circulations on $\Gamma$ represent precisely the $1$-cycles of $\Gamma$  (viewed as $1$-complex). The classical result in homology theory describes $H_1(\Gamma, \MZ)$ as the abelianization of the fundamental group $\pi_1(\Gamma)$, which, in this case, is isomorphic to  the free group $N$.  Putting all these together one has another geometric proof of the Fox theorem, as well as the description of the kernel of the Magnus embedding.

In Section  \ref{subsec:geodesics}  we describe  geodesics in groups $F/N^\prime$ as Euler tours in some finite subgraphs of $\Gamma$ generated by the supports of the flows  of the elements of $F/N^\prime$ on $\Gamma$. The description is geometric, explicit,  and it gives a natural way to compute the geodesic length of elements. In this part geometric ideas seem unavoidable.  However, this simplicity becomes treacherous when one concerns the efficiency of computations.

 We prove  that BGLP (relative to the standard basis)   is NP-complete even in $S_{r,2}$.
  Consequently,  the problems GP and GLP are NP-hard in $S_{r,2}$. To show this we construct a polynomial-time reduction of the Rectilinear Steiner Tree Problem (RSTP), which is NP-complete,  to BGLP in $S_{r,2}$.  The necessary information on RSTP is outlined in Section  \ref{subsec:RSTP} and   the proof of the main theorem is in Section \ref{subsec:MT}. Notice, that in \cite{DLS} GLP was claimed to be polynomial time decidable in arbitrary finitely generated free solvable groups, but the  argument turned out to be fallacious.

In the second half of the 20th century free solvable groups, as well as, solvable wreath products of groups and finitely generated metabelian groups, were intensely studied, but mostly from the view-point of combinatorial groups theory.
Now they stand at the heart of research  in various areas of algebra. On the one hand, the rejuvenated interest to these groups stems from   random walks on groups and,
cohomology theory.  For example, wreath products of abelian groups give
 exciting examples and  counterexamples to several conjectures on the numerical
characteristics of random walks. It seems, the main reasons that facilitate research here come from some  paradoxical properties of the groups itself: all
these groups are amenable (as solvable group), but they have exponential growth and may  have
nontrivial Poisson boundary  \cite{Vershik_Kaimanovich:1983}, etc. These groups, contrary to, say, free nilpotent groups, may have irreducible unitary representations with nontrivial cohomology.   Some  numerical characteristics of these groups are very
intriguing, giving  new exciting examples in the quantitative group theory.
For example,  metabelian "lamplighter" groups  have intermediate growth of the drift, positive entropy,  etc. These groups were  intensively studied recently (see papers \cite{Vershik_Kaimanovich:1983,Ershler3:2003} and  the  bibliography in the latter).

On the other hand, metabelian groups are currently at the focus of a very active research in geometric groups theory.   In  1983 Gromov proposed a program for studying finitely generated groups as geometric objects \cite{Gromov_inf_groups__geom_objects:1983}. One of the principal directions of this program  is the classification of finitely generated groups up to quasi-isometry. It follows from Gromov's result on groups with polynomial growth
\cite{Gromov_pgrowth:1981} that a  group quasi-isometric to a nilpotent group is virtually
nilpotent.  In the case of solvable groups the situation is much less known.  Erschler   shown in \cite{Dyubina:2000} that a group quasi-isometric to a solvable group may
be  not virtually solvable. Thus, the class of virtually solvable
groups is not closed under quasi-isometry. On the other hand there are interesting classes of solvable non-polycyclic groups that are
quasi-isometrically rigid, for example, solvable Baumslag-Solitar groups (Farb and Mosher \cite{Farb_Mosher1:1998,Farb_Mosher2:1999}). We refer to the papers \cite{Mosher_Sageev_Whyte:2003}  and \cite{Eskin_Fisher_Whyte:2007} for some recent results in this area.

 It seems timely to try extend the results of this paper to the  classes of solvable groups mentioned above.
There are many interesting open questions concerning computational complexity of algorithmic problems in these classes of  solvable groups, we discuss some of them in Section \ref{se:open-problems}.

All polynomial time algorithms presented in this work are
implemented and available at \cite{CRAG}.

\medskip

\noindent {\bf Acknowledgement}
\medskip

We would like to thank M.Sapir and B.Steinberg who brought to our attention some relevant geometric ideas from semigroup theory.

\section{Preliminaries}
  \label{se:preliminary}

\subsection{The Word Problem}
  \label{subse:WP}

Let $F = F_r =  F(X)$ be a free group with a basis  $X = \{x_1, \ldots,
x_r\}$. A subset $R \subseteq F$ defines a {\em presentation} $P = \langle X \mid R \rangle$ of a group
 $G = F / N$ where $N = ncl(R)$ is the normal closure of $R$ in $F$. If $R$ is finite (recursively enumerable) then the presentation is called finite (recursively enumerable).

The Word Problem $WP$ for $P$ is termed {\em decidable}
if the normal closure $N$ is a decidable subset of $F(X)$, i.e.,
there exists  an algorithm $\CA$ to decide whether a given word $w
\in F(X)$ belongs to $N$  or not. The time function
$T_\CA: F(X) \rightarrow \mathbb{N}$  of the algorithm $\CA$ is
defined  as the number of steps required for $\CA$ to halt on an
input $w \in F(X)$. We say that the Word Problem for $P$  is
decidable in polynomial  time if there exists a decision algorithm
$\CA$, as above, and constants $c, k \in \mathbb{N}$ such that
 $$T_\CA(w) \leq c|w|^k$$
 for every $w \in F(X)$ (here $|w|$ is the length of the word $w$).
 In this case we say that the time complexity of WP for $P$ is $O(n^k)$.

\subsection{Free solvable groups and the Magnus embedding}
  \label{subse:free-solvable}  \label{subsec:Magnus-embedding}

For a free group $F = F(X)$ of rank $r$ denote by
$F^{(1)} = F' = [F,F]$ the {\em derived} subgroup  of $F$, and by
$F^{(d)} = [F^{(d-1)}, F^{(d-1)}]$ -  the {\em $d$-th derived subgroup}
of $F$, $d\geq 2$. The quotient group $A_r = F_r / F_r'$ is a
{\em free abelian group} of rank $r,$ $M_r = F_r / F_r^{(2)}$ is a {\em free
metabelian group} of rank $r,$ and $S_{r,d} = F_r / F_r^{(d)}$ is a
{\em free solvable group} of rank $r$ and class $d$. In the sequel we
usually identify the set $X$ with its canonical images in $A_r,
M_r$ and $S_{r,d}$.

One of the most powerful approaches  to study free solvable groups is via the, so-called,  Magnus embedding.
To explain we need to introduce some notation. Let
$G = F/N$ and  $\mathbb{Z}G$  the group ring of $G$ with integer coefficients. By  $\mu:F \rightarrow G$ we denote the canonical factorization epimorphism, as well its linear extension to $\mu: \mathbb{Z}F \rightarrow \mathbb{Z}G$. Let $T$ be a free (left) $\mathbb{Z}G$-module of rank $r$ with a basis $\{t_1, \ldots, t_r\}$. Then the set of matrices
$$M(G) = \left(
 \begin{array}{ll}
  G & T\\
  0 & 1
  \end{array}
  \right)
= \left \{ \left(
 \begin{array}{ll}
  g & t\\
  0 & 1
  \end{array}
  \right) \mid  g \in G, t \in T \right \}
  $$
forms a group with respect to the matrix multiplication. It is easy to see that the group $M(G)$ is a discrete wreath product $M(G) = A_r wr G$ of the free abelian group $A_r$ and $G$.

In \cite{Magnus:1939} Magnus showed that the homomorphism $\phi:F \rightarrow M(G)$ defined by
$$ x_i \rightarrow \left(
 \begin{array}{ll}
  x_i^\mu & t_i\\
  0 & 1
  \end{array}
  \right),
  \ \ \ i = 1, \ldots, r,
  $$
satisfies $\ker \phi = N^\prime$. It follows that $\phi$ induces a monomorphism
    $$\psi: F/N^\prime \hookrightarrow M(F/N),$$
which is now called the {\em Magnus embedding}.

The Magnus embedding allows one to solve WP in the group
$F/N^\prime$ if WP in $G = F/N$ is decidable. Indeed,
given a word $w \in F(X)$ one can compute its image
$\phi(w)$ in $M(G)$ (multiplying the  images of
the letters in $w$) and then, using
a  decision algorithm for WP in $G$, check if the resulting
matrix $\phi(w)$ is the identity matrix or not.  To estimate the
complexity of such an algorithm notice first, that  the
coefficients from $\mathbb{Z}G$ that occur in the upper-right
corner of the matrix $\phi(w)$ have  $O(|w|)$ summands.
Secondly, to check whether or not  an element $h = m_1v_1 + \ldots + m_kv_k \in \mathbb{Z}G$,
where $m_i\in \MZ$ and $v_i \in G$ are given as words in the generators $X$ from $G$,
is trivial in $\mathbb{Z}G$  it requires about $O(k^2)$
comparisons of the type $v_i = v_j?$ in $G$. This gives an estimate
for the time function $T^\prime$ of WP in $F/N^\prime$ via the time function $T$ for WP in $F/N$:
    $$T^\prime(n) = O(rn^2T(n)),$$
where $n = |w|$.
Since WP in $A_r$ can be decided in linear time the estimate above shows that the complexity of WP in $M_r$ is  $O(rn^3)$. Moreover, induction on the solvability class $d$ gives a polynomial  estimate  $O(r^{d-1}n^{2d-1})$ for WP in the free solvable group $S_{r,d}$. Thus, the Magnus embedding gives a straightforward polynomial time (in $r$ and $n$) decision algorithm for $WP$ in $S_{r,d}$, but the degree of the polynomial grows with $d$. In particular, this algorithm is not polynomial as a uniform algorithm on the whole class of free solvable groups.

\subsection{Free Fox derivatives}
 \label{subsec:Fox-der}

Let $F = F_r(X)$ be a free group of rank $r$ with a basis $X = \{x_1,
... , x_r\}$.
The trivial group  homomorphism $F \rightarrow 1$ extends to a ring
homomorphism $\varepsilon : \mathbb{Z}F \rightarrow \mathbb{Z}$. The
kernel of $\varepsilon$ is called {\em the fundamental ideal} $\Delta_F$
of $\mathbb{Z}F$, it is a free (left) $\mathbb{Z}F$-module freely generated
by elements $x_1 - 1,  \ldots, x_r - 1$.

In \cite{Fox_calc1,Fox_calc2,Fox_calc3,Fox_calc4} R.~Fox introduced and gave a thorough account of the free differential
calculus in the group ring $\mathbb{Z}F$. Here we recall some notions and results referring to books
\cite{Crowell_Fox, Birman_book,Gupta} for details.

A map $D: \mathbb{Z}F \rightarrow \mathbb{Z}F$ is called a
{\it derivation} if it satisfies the following conditions:
\begin{enumerate}
    \item[\bf(D1)]
$D(u + v) = D(u) + D(v)$;
    \item[\bf(D2)]
$D(uv) = D(u) v^\varepsilon + uD(v)$, where $\varepsilon$ is the
ring homomorphism defined above.
\end{enumerate}

For every $x_i \in X$ there is a unique derivation, the so-called, a free partial
derivative $\partial / \partial x_i$,  such that
  $ \partial x_j/ \partial x_i = \delta _{ij}$, where $\delta_{ij}$ is the Kronecker's delta.
 It turned out that for every $u \in \mathbb{Z}F$
\begin{equation} \label{eq:Fox}
u - u^{\varepsilon } =  \partial u/ \partial x_1 (x_1 - 1) + ... +
\partial u/\partial x_r (x_r - 1 ).\end{equation}
 Since $\Delta_F$ is a free $\mathbb{Z}F$-module the equality
 (\ref{eq:Fox})
 gives another definition of the partial derivatives.

Condition (D2) implies the following useful formulas, that allow
one to compute easily partial derivatives of elements of
$\mathbb{Z}F$
\begin{equation} \label{eq:formula-inverse}
 \partial x_j^{-1}/\partial x_i = -\delta_{ij}x_j^{-1}.
\end{equation}
and, hence, for a word $w = x_{i_1}^{\varepsilon_1} \ldots
x_{i_n}^{\varepsilon_n} \in F(X)$ one has
\begin{equation} \label{eq:derivative}
\begin{array}{cc}
     \partial w/\partial x_i = \Sigma_{j = 1}^n x_{i_1}^{\varepsilon_1} \ldots x_{i_{j-1}}^{\varepsilon_{j-1}}(\partial x_{i_j}^{\varepsilon_j}/ \partial x_i) = \\
\end{array}
\end{equation}
     $$= \sum_{1\le j\le n, ~i_j = i,~ \varepsilon_j=1} x_{i_1}^{\varepsilon_1} \ldots x_{i_{j-1}}^{\varepsilon_{j-1}} - \sum_{1\le j\le n, ~ i_j = i,~ \varepsilon_j=-1} x_{i_1}^{\varepsilon_1} \ldots x_{i_{j}}^{\varepsilon_{j}}$$

The following result is one of the principle technical tools in this area, it follows easily from the Magnus embedding theorem, but in the current form it is due to Fox \cite{Fox_calc1,Fox_calc2,Fox_calc3,Fox_calc4}.

\medskip
\noindent{\bf Theorem} [Fox] {\it Let $N$ be a normal subgroup of $F$ and $\mu:F \to F/N$ the canonical epimorphism.   Then for every $u \in F$ the following equivalence holds:
$$\forall i  \ \left (\partial u/\partial x_i \right )^{\mu} = 0 \Longleftrightarrow u \in [N,N].$$
}

\medskip

In particular, for $N = F^{(d)}$ the standard epimorphism $\mu:F \rightarrow S_{d} = F/F^{(d)}$ gives rise to a ring homomorphism
 $\mu:\mathbf{Z}F  \rightarrow \mathbf{Z}S_{d}$ such that
  \begin{equation}
 \label{eq:main-Fox}
 F^{(d+1)} = \{ u \in F \mid  (\partial u/ \partial x_i)^{\mu} =
 0 \ for \ i = 1, \ldots,r\}.
 \end{equation}
 Composition of $\partial/\partial x_i$ with $\mu$ gives  an {\it induced}  partial derivative
  $\partial^\mu /\partial x_i: \mathbf{Z}F  \rightarrow \mathbf{Z}S_{d},$
 which we often denote again by $\partial /\partial x_i $ omitting $\mu$ (when it is clear from the context).

 Partial derivatives $\partial^\mu /\partial x_i$ are useful when computing images under the Magnus embedding. Indeed,
 by induction on the length of $w \in F$ it is easy to show that the image of $w$, under the Magnus embedding $\phi:F/N^\prime \rightarrow M(F/N)$, can be written as follows
 $$ w^\phi  =  \left(
 \begin{array}{ll}
  w^\mu &  \sum_{i = 1}^r \partial^\mu w/\partial x_i \cdot t_i\\
  0 & 1
  \end{array}
  \right).
  $$
This shows that the faithfulness of the Magnus embedding is, in fact, equivalent to the  Fox Theorem above.

\subsection{Flows on $F/N$}
\label{subsec:flows}

In this section we relate flow networks on the Cayley graph of $F/N$ to the elements of $F/N^\prime$.

Let $X=\{x_1,\dots , x_r\}$ be a finite alphabet.
An $X$-\emph{labeled directed graph} $\Gamma$ (or
$X$-\emph{digraph}) is a pair of sets $(V,E)$ where
the set $V$ is called the {\em vertex set} and the set
$E \subseteq V \times V \times X$ is called the {\em edge set}.
An element $e = (v_1,v_2,x) \in E$ designates an edge with
the origin $v_1$ (also denoted by $o(e)$), the terminus $v_2$
(also denoted by $t(e)$), labeled by $x$.
If for $e\in E$ we have $o(e)=t(e)$ then we say that $e$ is a $\emph{loop}$.
The graph $\Gamma$ can be finite or infinite.

\begin{example}
The Cayley graph $\Gamma(G,X)$ of the group $G = F/N$ is an $X$-digraph.
\end{example}

Given an $X$-digraph $\Gamma$, we can make $\Gamma$
into a directed  graph labeled by the alphabet $X^{\pm 1} =X\cup
X^{-1}$. Namely, for each edge $e = (v_1,v_2,x)$ of $\Gamma$ we introduce a
formal inverse $e^{-1} = (v_2,v_1,x^{-1})$.
For the new edges $e^{-1}$ we set $(e^{-1})^{-1}=e$.
The new graph, endowed with this additional structure, is
denoted by $\hat \Gamma$. In fact in many instances we
abuse notation by disregarding the difference between $\Gamma$ and
$\hat \Gamma$.

\begin{remark}
If $X$ is a generating set of $G$ such that $X \cap X^{-1} = \empty$ then   $\hat \Gamma(G,X)$ is the Cayley graph $\Gamma(G,X^{\pm 1})$ of $G$ relative to the generating set $X^{\pm 1}$.
\end{remark}

The edges of $\hat \Gamma$ inherited from $\Gamma$
are called \emph{positively oriented} or \emph{positive}.
The formal inverses of positive edges in $\hat \Gamma$ are called
 \emph{negatively oriented} or \emph{negative}.
The edge set of $\hat \Gamma$ splits in a disjoint union
$E(\hat\Gamma)=E^+(\Gamma) \sqcup E^-(\Gamma)$ of the sets of
positive and negative edges.

The use of $\hat \Gamma$ allows us to define the notion of a
\emph{path} in $\Gamma$. Namely, a \emph{path} $p$ in $\Gamma$ is
a sequence of edges $p = e_1,\dots , e_k$ where each $e_i$ is an
edge of $\hat \Gamma$ and the origin of each $e_i$ (for $i>1$) is the
terminus of $e_{i-1}$. In this situation we say that the
\emph{origin} $o(p)$ of $p$ is $o(e_1)$ and the \emph{terminus}
$t(p)$ is $t(e_k)$. The \emph{length} $|p|$ of this path is
set to be $k$.
Also, such a path $p$ has a naturally defined label
$\nu(p)=\nu(e_1)\dots \nu(e_k)$. Thus $\nu(p)$ is a word in
the alphabet $\Sigma=X\cup X^{-1}$. Note that it is possible that
$\nu(p)$ contains subwords of the form $a a^{-1}$ or $a^{-1} a$ for
some $a\in X$. If $v$ is a vertex of $\Gamma$, we will consider the sequence
$p=v$ to be a path with $o(p)=t(p)=v$, $|p|=0$ and
$\nu(p)=1$ (the empty word).

In general, one can consider labels in an arbitrary inverse semigroups, the construction above  applies to this case as well. In particular, we will consider directed graphs with labels in  $\MZ$.
We consider also digraphs with no labels at all (to unify terminology, we view them sometimes, as labeled in the trivial semigroup $\{1\}$), the construction above still applies.

Let $\Gamma = (V,E)$ be an $X$-digraph with two
distinguished vertices $s$ (called {\em source}) and $t$ (called {\em sink})
from $V$. Recall that a {\it flow} (more precisely
$\MZ$-flow) on $\Gamma$ is a function $f: E \rightarrow \MZ$ such that
\begin{itemize}
    \item[\bf(F)]
for all $v \in V - \{s,t\}$ the equality $\sum_{o(e) = v} f(e) - \sum_{t(e) = v} f(e) = 0$ holds.
\end{itemize}
The number $f^\ast(v) = \sum_{o(e) = v} f(e) - \sum_{t(e) = v}f(e)$
is called the {\it net flow} at $v \in V$. The condition (F) is
often referred to as the {\em Kirchhoff law} (see, for example,
\cite{Bollobas_intr_1990,Diestel}) or a {\em conservation law} \cite{Chartrand_Oellermann}.

For the digraph  $\hat \Gamma$ the definition above can be formulated in the following equivalent way, which is the standard one in flow networks:
\begin{itemize}
    \item[\bf(F1)]
$f(e) = - f(e^{-1})$ for any $e \in E$.
    \item[\bf(F2)]
$\sum_{o(e) = v} f(e) = 0$ for all $v \in V - \{s,t\}$.
\end{itemize}
Here the net flow at $v$ is equal to $f^\ast(v) = \sum_{o(e) = v} f(e)$.

Usually a flow network comes equipped with a  {\it capacity} function $c:E \rightarrow \MN$,  in which case a flow $f$ has to satisfy the {\it capacity restriction}
\begin{itemize}
    \item[\bf(F3)]
$f(e) \leq c(e)$ for all $e \in E$.
\end{itemize}
In the sequel we do not make much use of the capacity
function (it occurs in an obvious way), so in most cases
we consider flows on graphs $\Gamma$ satisfying  the Kirchhoff law  (F)
(or, equivalently, on graphs  $\hat \Gamma$ satisfying (F1) and (F2)).

A flow $f$ is called a {\it circulation} if (F) holds for all vertices from $V$ (including $s$ and $t$).

\begin{example}
Let $\Gamma = \Gamma(G,X)$ be the Cayley graph of $G = F/N$
relative to the generating set $X$. The constant function
$f: E(\Gamma) \rightarrow \{1\}$ defines a circulation on
$\Gamma$, since for every vertex $g \in V(\Gamma)$ and every
label $x \in X$ there is precisely one edge $(gx^{-1}, g)$
with label $x$ incoming into $g$ and precisely one edge $(g,gx)$
with the label $x$ leaving $g$.
 \end{example}

An important class of flows on $\Gamma = \Gamma(G,X)$
comes from paths in $\Gamma$. A  path $p$ in $\Gamma$
defines  an integer-valued  function
$\pi_p:E(\Gamma) \rightarrow \MZ$, such that for an edge $e$
$\pi(e)$ is the algebraic sum (with respect to the orientation)
of the number of times the path $p$ traverses $e$, i.e.,
each traversal of $e$ in the positive direction  adds $+1$,
and in the negative direction adds $-1$.  It is obvious that
$\pi_P$ is a flow on $\Gamma$ with the source $o(p)$
(the initial vertex of $p$) and the sink $t(p)$
(the terminal vertex of $p$). Notice that $\pi_p$ satisfies
also the following conditions:
\begin{itemize}
    \item[\bf(F4)]
either $\pi_p$ is a circulation (iff $p$ is a closed path), or
$f^\ast(s)  = 1, f^\ast(t) = -1$.
    \item[\bf(F5)]
$\pi_p$ has finite support, i.e, the set
$supp(\pi) = \{e \in E \mid \pi(e) \neq 0\}$ is finite.
\end{itemize}
We say that a flow $\pi$ on $\Gamma$ is {\it geometric} if it satisfies conditions (F4) and (F5).

It is easy to see that the set $\CC(\Gamma)$ of all circulations on $\Gamma$ forms  an abelian  group with respect to the operations (here $f,g \in \CC(\Gamma)$):
\begin{itemize}
    \item
$(f+g)(e) = f(e) + g(e)$,
    \item
$(-f)(e) = - f(e)$.
\end{itemize}
Meanwhile, the  set $\mathcal{GC}(\Gamma)$ of all geometric circulations is a subgroup of $\CC(\Gamma)$.
On the other hand,  the sum $f + g$ of two geometric flows gives a geometric flow only if the sink of $f$ is equal to the source  of $g$, or either $f$ or $g$ (or both) is a circulation. In fact, the set $\mathcal{GF}(\Gamma)$  of all geometric flows is a groupoid.

Let $\Pi(\Gamma)$ be the fundamental groupoid of paths in $\Gamma$. Then the map   $\sigma: \Pi(\Gamma) \rightarrow \mathcal{GF}(\Gamma)$ defined for $p \in \Pi(\Gamma)$ by $\sigma(p) =  \pi_p$ is a morphism in the category of  groupoids, i.e., the following holds (here $p, q \in \Pi(\Gamma)$):

\begin{itemize}
\item $\pi_{pq} = \pi_p + \pi_q$, if  $pq$ is defined in $\Pi(\Gamma)$,
\item $\pi_{p^{-1}} = -\pi_p$.
\end{itemize}

Now we will show that every geometric flow $\pi$ on $\Gamma$ can be realized as a path flow $\pi_p$ for a suitable path $p$.

\begin{lemma}
\label{le:geometric-flow}
Let $\pi$ be a geometric  flow on $\Gamma$. Then there exists a path $p$ in $\Gamma$ such that $\pi = \pi_p$.
\end{lemma}
\begin{proof}
Let $\pi$ be a geometric  flow on $\Gamma$ with the source $s$ and the sink $t$. Denote by $\Gamma_\pi$ the subgraph of $\Gamma$ generated by $supp(\pi) \cup \{s,t\}$. Suppose $Q$ is a subgraph of $\Gamma$ such that  $\Delta = \Gamma_\pi \cup Q$ is a connected graph (every two vertices are connected by a path in $\hat \Gamma$). Clearly, $\pi$ induces a flow on $\Delta$.  Now  we construct another $X$-digraph $\Delta^\ast$ by adding new edges to $\Delta$ in the following manner. For every edge $e \in E(\Delta)$ with $|\pi(e)| > 1$ we add extra $|\pi(e)| - 1$ new edges $e^{(1)}, \ldots, e^{|\pi(e)| - 1}$ from $o(e)$ to $t(e)$ if $\pi(e) > 0$ and from $t(e)$ to $o(e)$, if $\pi(e) < 0$. We label the new edges by the same label if $\pi(e)  > 0$, and by its inverse, otherwise. If $\pi(e) = 0$ then we add a new edge $e^{-1}$  from  $t(e)$ to $o(e)$ with the inverse label.  In the case $|\pi(e)| = 1$ we do not add any new edges.  Notice that every vertex in $V(\Delta^\ast) -\{s,t\}$  has even directed degree (the number of incoming edges is equal to the number of outgoing edges). There are two cases to consider.

Case 1). Suppose  $\pi$ is a circulation. Then every vertex in $\Delta^\ast$ has even directed degree. Therefore, the digraph $\Delta^\ast$ has an Euler tour $p^\ast$, i.e., a closed path at $s$ that traverses every edge in $\Delta^\ast$  precisely once.
 Let $\phi:\Delta^\ast \rightarrow \hat{\Delta}$ be the morphism of $X$-digraphs that maps all the new edges $e^{(1)}, \ldots, e^{|\pi(e)| - 1}$ to their original edge $e$. Clearly, the image $p = \phi(p^\ast)$ is a path in $\hat{\Delta}$ such that $\pi_p = \pi$.

Case 2). Suppose $\pi$ is not a circulation.  Let $q$ be a path in  $\Delta^\ast$ from $s$ to $t$.  Then $\pi^\prime = \pi - \pi_q$ is a circulation. Hence by Case 1) there exists a path $p$ in $\Gamma$  such that  $\pi^\prime = \pi_p$. Therefore, $\pi_{pq} = \pi_p + \pi_q = \pi^\prime + \pi_q = \pi$, as required.

\end{proof}

\subsection{Geometric interpretation of Fox  derivatives}
 \label{subsec:geom-interpret}

In this section we give a geometric interpretation of Fox derivatives.

Let  $G = F/N$, $\mu:F \rightarrow F/N$  the canonical epimorphism,
and  $\Gamma = \Gamma(G,X)$ the Cayley graph of $G$ with respect to
the generating set $X$.
A word  $w \in F(X)$ determines a unique path $p_w$ in $\Gamma$ labeled by $w$
which starts at 1 (the vertex corresponding to the identity of $G$).
As we mentioned in Section \ref{subsec:flows} the
path $p_w$ defines  a geometric flow $\pi_{p_w}$ on $\Gamma$ which we denote by $\pi_w$.

\begin{lemma}
\label{le:deriv-geom}
Let $w \in F = F(X)$. Then for any  $g\in F/N$ and  $x \in X$ the value of $\pi_w$
on the edge $e = (g,gx)$ is equal to the coefficient in front of $g$ in the Fox derivative
$(\partial w/\partial x)^\mu \in \MZ G$, i.e.,
    $$(\partial w/\partial x)^\mu  = \sum_{g \in G, x\in X}\pi_w(g,gx) g$$
\end{lemma}

\begin{proof}
Follows by induction on the length of $w$  from formulas (\ref{eq:derivative}).
\end{proof}

\begin{figure}[htbp]
\centering
\epsfig{file=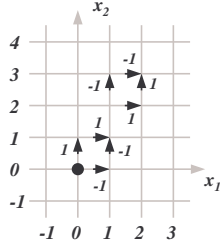, scale=0.6}
\caption{\label{fi:diagram1}
The values of $\pi_w$ for $w = x_2 x_1 x_2 x_1 x_2 x_1^{-1} x_2^{-3}
x_1^{-1}$ on $(x_1,x_2)$-grid. In this case  $\partial w/ \partial
x_1 = -1 + x_2 - x_1 x_2^3 + x_1 x_2^2$ and $\partial w/ \partial
x_2 = 1 - x_1 + x_2^2 x_2^2 - x_1 x_2^2$.}
\end{figure}

Figure \ref{fi:diagram1} is an example
for $F = F(\{x_1,x_2\})$ and $N = F^\prime$.
Non-zero values of $\pi_w$ are shown as weights on the edges (zero weights are omitted).

The following theorem has been proven in \cite{DLS,Vershik_Dobrynin:2004} using a  homological argument similar to the one in  Proposition \ref{pr:geometric-circ}. Here we give a short independent proof based on the Fox theorem.
\begin{theorem}
\label{th:pi}
 \cite{DLS,Vershik_Dobrynin:2004} Let $N$ be a normal subgroup of $F$ and $\sim_N$ an equivalence relation on $F$ defined by
  $$u \sim_N v  \Longleftrightarrow \pi_u = \pi_v.$$
  Then
$F/N^\prime  = F/ \sim_N$.
\end{theorem}

\begin{proof}
Let $u,v \in F$. Suppose $u = v$ in $F/N^\prime$.
Then $uv^{-1} \in N^\prime$ hence by Fox theorem $\partial^\mu (uv^{-1})/\partial x = 0$ for every $x \in X$
(here by $\partial^\mu/\partial x$ we denote the canonical
image of $\partial /\partial x$ in the group ring $\MZ(F/N)$).
Observe that
\begin{equation}\label{eq:product-inverse-derivatives}
    \frac{\partial }{\partial x}(uv^{-1}) = \frac{\partial u}{\partial x} - uv^{-1}\frac{\partial v}{\partial x}.
\end{equation}
Hence in $\MZ(F/N)$
    $$0 = \frac{\partial^\mu}{\partial x}(uv^{-1}) = \frac{\partial^\mu u}{\partial x} - \frac{\partial^\mu v}{\partial x},$$
so, by Lemma \ref{le:deriv-geom} $\pi_u = \pi_v$, as claimed.

To show the converse, notice first that $\pi_u = \pi_v$ implies
that $u = v$ in $F/N$. Indeed, it can be seen from the definition
of $\pi$ but also follows from (\ref{eq:Fox}) and Lemma \ref{le:deriv-geom}
since in this case
    $$(u^\mu - v^\mu) - (u-v)^\varepsilon = \sum_{x \in X} \frac{\partial^\mu}{\partial x}(u-v) \cdot (x-1) = 0$$
which implies $u^\mu = v^\mu$. Now by (\ref{eq:product-inverse-derivatives})
    $$\frac{\partial^\mu }{\partial x}(uv^{-1}) = \frac{\partial^\mu u}{\partial x} - \frac{\partial^\mu v}{\partial x} = 0,$$
and, hence, by Fox theorem $uv^{-1} \in N^\prime$.

 \end{proof}

\begin{remark}
 Theorem \ref{th:pi} relates the algebraic and geometric points of view on derivatives.
  One can prove this theorem (see Section \ref{subsec:homology}) using a pure topological argument, then the Fox theorem, as well as, the Magnus embedding, come along as easy corollaries.
\end{remark}

\subsection{Geometric circulations and the first homology group of $\Gamma$}
\label{subsec:homology}

 We describe here  geometric circulations on $\Gamma = \Gamma(G,X)$ in a pure topological manner. For all required notions and results on homology of simplicial complexes we refer to \cite{Hatcher} or \cite{SingerThorpe67}.

 In this section  we view $\Gamma$ as an infinite $1$-complex.

\begin{proposition}
\label{pr:geometric-circ}
Let $G = F/N$,  $\Gamma = \Gamma(G,X)$, and  $\sigma:   \pi_1(\Gamma) \rightarrow \mathcal{GC}(\Gamma)$ be a map defined by $\sigma(p) =  \pi_p$ for $p\in \pi_1(\Gamma)$ . Then:
 \begin{itemize}
 \item $\sigma$ is an epimorphism of groups;
 \item every geometric circulation  on $\Gamma$ defines a $1$-cycle on $\Gamma$;
 \item $H_1(\Gamma,\MZ) \simeq \mathcal{GC}(\Gamma)$;
 \item $\pi_1(\Gamma)  \simeq N$ and $\ker \sigma  = N^\prime$.
 \end{itemize}

\end{proposition}

\begin{proof}
It was mentioned already in Section \ref{subsec:flows}  that the map   $p \rightarrow  \pi_p$ is a morphism from the  fundamental groupoid $\Pi(\Gamma)$ of paths in $\Gamma$ into the groupoid of geometric flows  $\mathcal{GF}(\Gamma)$. Hence,  the restriction of this map onto the fundamental group $\pi_1(\Gamma)$ of $\Gamma$ gives a homomorphism of groups.
We have seen in Lemma \ref{le:geometric-flow} that $\sigma$ is onto.  This proves the first statement.

To see 2) observe first that a geometric circulation $f: E(\Gamma) \rightarrow \MZ$, viewed as a formal sum $\sum_{e \in E(\Gamma)}f(e)e$,  gives precisely a 1-chain in $\Gamma$ (see, for example, \cite{Hatcher}). Moreover,by definition, the net flow $f^\ast(v)$ at the vertex $v \in V(\Gamma)$ is the coefficient in front of $v$ in the boundary $\partial f$ of $f$. Therefore, $\partial f = 0$, so $f$ is a $1$-cycle.

3) follows easily from the 2). Indeed, $\Gamma$ is an $1$-complex, so there are no non-trivial $1$-boundaries in $\Gamma$. In this event $H_1(\Gamma, \MZ)$ is isomorphic to the group $\mathcal{GC}(\Gamma)$ of $1$-cycles, as claimed.

4) It is a classical  result that the kernel of $\sigma: \pi_1(\Gamma) \rightarrow H_1(\Gamma,\MZ)$ is equal to the derived subgroup of $\pi_1(\Gamma)$ (see \cite{Hatcher}).
To prove 4) it suffices to notice that $\pi_1(\Gamma) \simeq N$, which is easy.
\end{proof}

\begin{remark}
 Proposition \ref{pr:geometric-circ} gives a simple geometric proof of Theorem \ref{th:pi}, that relates the algebraic and geometric points of view on derivatives.
  Now one can derive the Fox theorem from the geometric argument above and then obtain the description of the kernel of the Magnus embedding, as an easy corollary.
\end{remark}

\subsection{Geodesics in $F/N^\prime$}
\label{subsec:geodesics}

Let $G = F/N$ and $\mu:F(X) \rightarrow G$ the canonical epimorphism.
In this section  we describe geodesics of elements of the
group $H = F/N^\prime$ relative to the set of generators $X^\mu$.

It is convenient to view the free group $F = F(X)$ as the
set of all freely reduced words in the alphabet
$X^{\pm 1} = X \cup X^{-1}$ with the obvious multiplication.

To describe geodesics in $H$ (relative to $X$) of a given
word $w \in F$  we need a construction from Lemma \ref{le:geometric-flow}.
Recall that $p_w$ is the path in the Cayley graph $\Gamma = \Gamma(G,X)$
from $1$ to $w^\mu$ with the label $w$ and $\pi_w$ the induced
geometric flow on $\Gamma$ with the source $1$ and the
sink $w^\mu$. By $\Gamma_w$ we denote the subgraph of $\Gamma$
generated by $supp(\pi_w) \cup \{1, w^\mu\}$. Suppose $Q$ is a
finite subgraph of $\Gamma$ such that  $\Delta = \Gamma_w \cup Q$
is a connected graph and $Q$ has the least number of edges among
all such subgraphs.  It follows from minimality of $Q$ that every
connected component of $Q$ is a tree. Moreover, if in the graph
$\Delta = \Gamma_w \cup Q$ one collapses every connected component
 $\Gamma_w$ to a point, then the resulting graph is a tree. We
 refer to $Q$ as a {\it minimal forest} for $w$. In general, there
 could be several minimal forests  for $w$.

Similarly as in the proof of Lemma \ref{le:geometric-flow},
we construct a finite   $X$-digraph $\Delta^\ast$ by
replicating every edge $e \in E(\Delta)$ with  $||\pi_w(e)| - 1|$
new edges in such a way that every vertex in $V(\Delta^\ast) -\{1,w^\mu\}$
has even directed degree  and the map, that sends every
replica of an edge $e$ back to $e$ (or $e^{-1}$ depending on
the orientation) is  a morphism of $X$-labeled digraphs
$\phi:\Delta^\ast \rightarrow \hat{\Delta}$. There are two cases.

\begin{itemize}
    \item
{\sc Case I.} Suppose that $p_w$ is a closed path in $\Gamma$, i.e., $w \in N$. In this case every vertex in $\Delta^\ast$ has even degree, so $\Delta^\ast$ has a Euler tour, say $p^\ast_{Q}$ at $1$. Denote by $p_Q$ the image $\phi(p^\ast_Q)$ of $p^\ast_{Q}$  under $\phi$. It follows from the construction (see Lemma \ref{le:geometric-flow}) that $p_Q$ is a closed path at $1$ in $\Gamma$ such that $\pi_w  = \pi_{p_Q}$. Therefore, if $w_{p_Q} \in F$ is the label of $p_Q$ then by Theorem \ref{th:pi} $w = w_{p_Q}$ in $H$. Moreover, since $p_Q$ is a Euler tour in $\Delta^\ast$ its length, hence the length of $w_{p_Q}$  is equal to
\begin{equation}\label{eq:path_length}
    |p_Q| = \sum_{e \in supp(p_w)} |\pi_w(e)| + 2|E(Q)|.
\end{equation}

    \item
{\sc Case II.} Suppose that $p_w$ is a not a closed path in $\Gamma$, i.e., $w \not\in N$.
By induction on $|w|$
it is easy to show that the vertices $1$ and $w^\mu$ belong to the same connected component of $\Gamma_w$.   Again, there exists  an  Euler's tour $p_Q^\ast$ in the graph $\Delta^\ast$ which starts  at the source and ends  at the sink. Clearly,  $\pi_{p_Q}$ satisfies the equality (\ref{eq:path_length}). If $u$ is the label of the path $\pi_{p_Q}$ then $\pi_u  = \pi_{p_Q} = \pi_w$ and $u$ is a geodesic word for $w$.

\end{itemize}

Now, with the construction in place, we are ready to characterize geodesics in $H$ of elements from $N$.

\begin{theorem}
\label{th:geodesics}
Let $H = F/N^\prime$ and $w \in F$. Then the following holds:
\begin{itemize}
\item if $Q$ is a minimal forest for $w$ then $w_{p_Q}$ is a geodesic for $w$ and
 $$l_X(w) = \sum_{e \in supp(p_w)} \pi_w(e) + 2|E(Q)|.$$
 \item every geodesic word for $w$ is equal (in $F$) to a word $w_{p_Q}$ for a suitable minimal forest $Q$ and an Euler path $p^\ast_Q$.
\end{itemize}
\end{theorem}
\begin{proof}
Let $u \in F$ be a geodesic word for $w^\mu$ in $H$. Observe, that  $\Delta = supp(\pi_u) \cup p_u$ is a connected subgraph of $\Gamma$ and
 $$|p_u| \geq \sum_{e \in supp(p_u)} \pi_u(e) + 2|E(\Delta - supp(\pi_u))|.$$
Now, by Theorem \ref{th:pi}  the equality $u = w$ in $H$ implies $\pi_u = \pi_w$. In particular, $supp(\pi_u) = supp(\pi_w)$. Hence
\begin{equation}
\label{eq:minimal-forest}
|p_u| \geq \sum_{e \in supp(p_w)} \pi_w(e) + 2|E(\Delta - supp(\pi_w))|.
\end{equation}
Since $u$ is geodesic for $w$ the number $|\Delta - supp(\pi_w)|$ is minimal possible, so
$Q = \Delta - supp(\pi_w)$ is a minimal forest for $w$. In fact, the equation \ref{eq:minimal-forest} shows that the converse is also true. This proves the theorem.
\end{proof}

The discussion above shows that GP is easy in $F/N^\prime$ provided one can solve the following problem efficiently.

\medskip
\noindent
 {\bf Minimal Forest Problem (MFP)}  Given a finite set of connected  finite subgraphs $\Gamma_1, \ldots, \Gamma_s$ in $\Gamma$ find a finite subgraph $Q$ of $\Gamma$ such that $\Gamma_1 \cup  \ldots \cup \Gamma_s \cup Q$ is connected and $Q$ has a minimal possible number of edges.

\medskip
\begin{proposition}
\label{pr:reduction-to-MFP}
GP in $F/N^\prime$ (relative to $X$) is linear time reducible to  MFP for $\Gamma(F/N,X)$.
\end{proposition}
 \begin{proof}
 Follows from the discussion above. Indeed, given a word $w \in F$ one can in linear time compute the flow  $\pi_w$ and find the connected components $\Gamma_1, \ldots, \Gamma_s$ of $supp(\pi_w) \cup\{o(p_w), t(p_w)\}$ in $\Gamma = \Gamma(F/N,X)$.  Then, solving MFP for these components, one gets the subgraph $Q$ which makes the graph  $\Gamma_1 \cup  \ldots \cup \Gamma_s \cup Q$ connected. Obviously, it takes linear  time to find a Euler path in the graph $\Delta$, hence to find a geodesic for $w$.
 \end{proof}

\section{The Word Problem in free solvable groups}
  \label{sec:WP-sovable-groups}

In this section we present fast algorithms to compute Fox derivatives of elements of  a free group $F$ in the group ring $\MZ S_{r,d}$ of a free solvable group $S_{r,d}$. As an immediate application,  we obtain a decision algorithm for WP in a free metabelian group $M_r$  with time complexity $O(rn\log_2n)$ and a decision algorithm for WP in $S_{r,d}, d \geq 3,$ with time complexity $O(rdn^3)$.  These are significant improvements in comparison with the known decision algorithms discussed in Introduction and Section \ref{subsec:Magnus-embedding}. As another application we get a fast algorithm to compute images of elements from $S_{r,d}$ under the Magnus embedding, which opens up an opportunity to use efficiently the classical techniques developed for wreath products.

\subsection{The Word Problem in free metabelian groups}
\label{sec:WP-metabelian}

In this section we compute Fox derivatives of elements of $F$ in the group ring $\MZ (F/F^\prime)$. Then we apply this  to WP in free metabelian groups.

Let $X = \{x_1, \ldots, x_r\}$, $F = F_r = F(X)$, $M  = M_r = F/F^{(2)}$, $A = A_r  = F/F^\prime$,  and $\mu:F \rightarrow A$ the canonical epimorphism. All Fox derivatives in this section are computed in the ring $\MZ A$.

Let $w \in F$. Then
$$
 \frac{\partial^\mu w}{\partial x_i} = \sum_{a \in A} m_{a,i}a, \ \ \ m_{a,i} \in \MZ.
$$
One can encode all the derivatives $\partial^\mu w /\partial x_i$ in one mapping $M_w: A \times \{1, \ldots,r\} \rightarrow \MZ$, where $M_w(a,i) = m_{a,i}$.  Let
$$
supp(M_w) = \{ (a,i) \mid M_w(a,i) \neq 0\},
$$
and $S_w$ the restriction of $M_w$ onto $supp(M_w)$. To compute Fox derivatives of $w$ we construct a sequence of finite maps $S_0 = \emptyset, S_1, \ldots, S_n = S_w$, as we read $w$.  On each step $k$ we either extend the domain $Dom(S_k)$ of $S_k$ or change the value of $S_k$ on some element from $Dom(S_k)$. To do this we need a data structure which allows one to do the following operations efficiently:

\begin{itemize}
    \item
for a given $(a,i)$  determine if $(a,i) \in Dom(S_k)$  or not;
    \item
add $(a,i)$ to $Dom(S_k)$ if $(a,i) \not \in Dom(S_k)$ and define $S_k(a,i) = q$ for some $q \in \MZ$;
    \item
 change the value of $S_k$ on $(a,i)$ if $(a,i) \in Dom(S_k)$.
\end{itemize}
Every element $a \in A$ can be written uniquely in the form
$a = x_1^{\delta_1(a)} \ldots x_r^{\delta_r(a)}$, where $\delta_i(a) \in \MZ$, so one may use the tuple of coordinates $\delta(a) = (\delta_1(a), \ldots, \delta_r(a))$ to represents $a$. It follows from the formula (\ref{eq:derivative}) for partial derivatives that for every $(a,i) \in Dom(S_w)$ the components $\delta_j(a)$ of $\delta(a)$ satisfy the inequality $|\delta_j(a)| \leq |w|$, as well as the values of $S_w$ and, hence, $|S_w(a,i)| \leq |w|$. Therefore, it takes $\lceil\log_2(|w|+1)\rceil$ bits to encode one coordinate $\delta_j(a)$ (one extra bit to encode the  sign $+$ or $-$), $r\lceil\log_2(|w|+1)\rceil$ bits to encode $\delta(a)$, and at most  $r\lceil\log_2(|w|+1)\rceil + \lceil\log_2r\rceil$ bits to encode $(a,i)$. We denote the binary word encoding $(a,i)$ by $(a,i)^\ast$.

\begin{figure}[htbp]
\centering
\epsfig{file=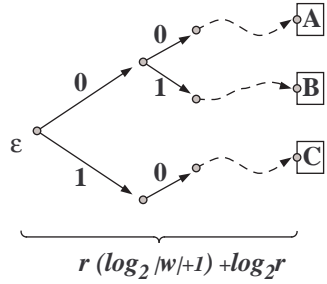, scale=1}
\caption{\label{fi:tree} The tree
$T_k$ representing the function $S_k$. }
\end{figure}

 Thus every  function  $S_k$ can be uniquely
represented by a directed $\{0,1\}$-labeled binary tree $T_k$ with the  root $\varepsilon$ and   with leaves labeled by integers such that
$(a,i)  \in Dom(S_k)$ if and only if there exists a path in $T_k$ from the root $\varepsilon$  to a leaf  labeled by the code of $(a,i)$ and such that the leaf of this path is labeled precisely by the integer $S_k(a,i)$.  Notice, that the  height of $T_k$ is equal to $r\lceil\log_2(|w|+1)\rceil+\lceil\log_2r\rceil$.
Such a tree is visualized schematically
in Figure \ref{fi:tree}.

\begin{remark}
\label{rem:data-structure}
It is clear that one can perform every operation mentioned  above on this data structure in at most
$r\lceil\log_2(|w|+1)\rceil+\lceil\log_2r\rceil$ elementary
steps.
\end{remark}

We use this data structure to
design  the following  algorithm for computing $S_w$.

\begin{algorithm}\label{al:metab-Fox}{\bf(Computing Fox derivatives in $F/F^\prime$)}
    \\{\sc Input.}
$r \in \MN$ and $w = x_{i_1}^{\varepsilon_1} \ldots x_{i_n}^{\varepsilon_n} \in F(X)$,
where $i_j \in \{1, \ldots, r\}$
and $\varepsilon_j = \pm 1$.
    \\{\sc Output.}
$S_w$.
    \\{\sc Computations.}
\begin{itemize}
    \item[A.]
Set $S = \emptyset$ and $\delta(a) = (0, \ldots, 0) \in \MZ^r$.
    \item[B.]
For $j = 1 ,\ldots, n$ do:
\begin{itemize}

    \item[(1)]
if $\varepsilon_j = 1$ then
\begin{itemize}
    \item
check if there is  a path (from the root to a leaf) labeled by $(\delta(a),i_j)^\ast$ in $S$;
    \item
if such a path does not exist in $S$ then create it,  add  it to $S$, and put $1$ as  the corresponding  value at the  new leaf;
    \item
if such a path exists in $S$ then add $1$ to  the value at its leaf.
\end{itemize}
    \item[(2)]
if $\varepsilon_j = -1$ then
\begin{itemize}
\item add $\varepsilon_{j}$ to  the $i_j$th coordinate of $\delta(a)$;
    \item
check if there is  a path (from the root to a leaf) labeled by $(\delta(a),i_j)^\ast$ in $S$;
    \item
if such a path does not exist in $S$  then create it, add it to $S$, and put $-1$ at the corresponding leaf;
    \item
if such a path exists in $S$ then subtract $1$ from the value at its leaf.
 \end{itemize}
\end{itemize}
\item[C.] Output $S_n$.
\end{itemize}
\end{algorithm}

\begin{theorem}
\label{th:comp-derivatives}
Given $r \in \MN$ and $w \in F$ Algorithm \ref{al:metab-Fox} computes all partial derivatives of $w$ in $\MZ A$ (the mapping  $S_w$)  in time $O(r |w| \log_2 |w|)$.
\end{theorem}
\begin{proof}
Using the formula (\ref{eq:derivative}) for partial derivatives it is easy to check that  given $w \in F$  Algorithm \ref{al:metab-Fox},  indeed, computes the mapping $S_w$.
 To  verify the complexity estimates observe, first that there are $|w|$ iterations in the part  $[B.]$ of Algorithm \ref{al:metab-Fox}. Each such iteration requires $O(r\lceil\log_2(|w|+1)\rceil)$ elementary steps  (see Remark \ref{rem:data-structure}),  so altogether one has $O(r|w|log_2|w|)$ as the time complexity estimate for Algorithm \ref{al:metab-Fox}, as claimed.

\end{proof}

\begin{algorithm}\label{al:wp_metab}{\bf(Word Problem in  free metabelian groups)}
    \\{\sc Input.}
$r \in \MN$ and $w \in F$.
    \\{\sc Output.}
$True$ if $w$ represents the identity in $M_r$ and $False$
otherwise.
    \\{\sc Computations.}
\begin{itemize}

     \item Apply Algorithm \ref{al:metab-Fox} to compute $S_w$.

     \item Check, looking at the values assigned to leaves of $S_w$, if all Fox derivatives $\partial w/\partial x_i$, $i = 1, \ldots,r$ are equal to $0$ or not.

     \item If  all the derivatives  are $0$ then output $True$. If there is a  non-zero derivative  then output $False$.
\end{itemize}
\end{algorithm}

\begin{theorem}
Algorithm \ref{al:wp_metab} solves the Word problem in a free
metabelian group $M_r$ in time $O(r |w| \log_2 |w|)$.
\end{theorem}

\begin{proof}
Follows from Theorem \ref{th:comp-derivatives} and the Fox Theorem (see Section \ref{subsec:Fox-der}).
\end{proof}

\subsection{The Word problem in free solvable groups}
 \label{subsec:free-solvable}

In this section we present an algorithm to compute all Fox derivatives of a word $w \in F$ in the group ring $\MZ S_{r,d-1}$, $d \geq 2$ in time $O(rd |w|^3)$. This gives a decision algorithm for WP in $S_{r,d}$ within time complexity $O(rd |n|^3)$.

Let $X = \{x_1, \ldots, x_r\}$, $F = F_r = F(X)$, $S = S_{r,d}  = F/F^{(d)}$, $d \geq 3$,   and $\mu:F \rightarrow S_{r,d-1}$ the canonical epimorphism. All Fox derivatives in this section are computed in the ring $\MZ S_{r,d-1}$.

In Section \ref{sec:WP-metabelian} we used a unique representation
of  elements $a \in A_r$ by their coordinate vectors $\delta(a)$ to compute  Fox derivatives in nearly linear time.  Now we do not have such normal forms of elements of $S_{r,d}$,  $d \geq 2$,  so our computations are slightly different, however, the general strategy is quite similar.  To speed up computations we use some data structures based on efficient partitioning techniques.

In general, let $G$ be a group generated by $X$   and $D$ a finite subset of $F(X)$. A {\em
$G$-partition} of $D$ is a partition of $D$ into a union of disjoint non-empty
subsets $D_i$ such that for any
$u,v \in D$, $u = v$ in $G$ if and
only if they belong to some subset $D_i$. Clearly, the $G$-partition of $D$ is unique. Observe that
if a group $H$ is a quotient  of $G$ then the
$G$-partition of  $D$ is the same or {\em finer} than the $H$-partition of $D$.

If the set $D$ is ordered, say $D = \{w_0, \ldots, w_n\}$ then the $G$-partition of $D$ can be represented by a function
$P:\{0,\ldots,n\} \rightarrow \{0,\ldots,n\}$ where $P(j) = i$ if and only if $w_i = w_j$ in $G$ and $i$ is the smallest with such property. Given the $G$-partition $P$ of $D$ one can arrange his data in such a way that it takes linear time (in the size of $j$) to compute   $P(j)$. In particular, given $i, j$ one can check in linear time if $w_i = w_j$ in $G$.    Also,  for a given word $w \in D$ one can find in linear time an index $i$ such that $w = w_i$.  These are the main two subroutines concerning partitions of $D$.

Let $w = x_{i_1}^{\varepsilon_1}
\ldots x_{i_k}^{\varepsilon_k}$, where $i_j \in \{1, \ldots, r\}$
and $\varepsilon_j = \pm 1$.
Put
\begin{equation}\label{eq:monomials}
    D_w = \{\varepsilon, ~~x_{i_1}^{\varepsilon_1}, ~~x_{i_1}^{\varepsilon_1} x_{i_2}^{\varepsilon_2}, ~~\ldots, ~~x_{i_1}^{\varepsilon_1} \ldots x_{i_n}^{\varepsilon_n}\} \subset F(X_r).
\end{equation}
We order the set $D_w$ as  follows $w_0 = \varepsilon, \ldots, w_n = w$.
Now, to check whether or not the derivative
$\partial w /\partial x_i$ is trivial in $\MZ S_{r,d-1}$ one has to determine which
pairs $(w_i,w_j)$ of elements from $D_w$ represent the same element in $S_{r,d-1}$ and then cancel out $w_i$ with $w_j$ in $\partial w /\partial x_i$ if they have the  opposite signs.

The goal  of Algorithm \ref{al:wp_solvable} below is to compute
$S_{r,d}$-partition for $D_w$.  This is performed in a sequence of
iterations. The algorithm starts out by computing
the $A_r$-partition of $D_w$. On the second iteration the algorithm
computes the $M_r$-partition of $D_w$. On the third step it
computes the  $S_{r,3}$-partition of $D_w$. It continues this way until it computes the $S_{r,d}$-partition of $D_w$.

To explain how the algorithm works  assume that the $S_{r,d-1}$-partition
of $D_w$ is given by the  partition function $P_{d-1}$ described
above. Notice, that the $S_{r,d}$-partition
of $D_w$ is the same or finer than the $S_{r,d-1}$-partition of $D_w$, since $S_{r,d-1}$
is a quotient of $S_{r,d}$. This shows that to construct  $S_{r,d}$-partition $P_d$ of $D_w$ one has only to compare elements from $D_w$ which are equal in $S_{r,d-1}$.
Suppose that $w_s, w_t \in D_w$, $s < t$ and $w_s = w_t$ in $S_{r,d-1}$. To check if $w_s = w_t$ in $S_{r,d}$  we compare all their Fox derivatives in $\MZ S_{r,d-1}$, so for every $k = 1, \ldots, r$ we compute the following differences:
    $$\partial w_s/\partial x_k - \partial w_t/\partial x_k = $$
    $$\sum_{1\le j\le s, ~i_j = k,~ \varepsilon_j=1} x_{i_1}^{\varepsilon_1} \ldots x_{i_{j-1}}^{\varepsilon_{j-1}}
    - \sum_{1\le j\le s, ~ i_j = k,~ \varepsilon_j=-1} x_{i_1}^{\varepsilon_1} \ldots x_{i_{j}}^{\varepsilon_{j}} - $$
    $$ - \sum_{1\le j\le t, ~i_j = k,~ \varepsilon_j=1} x_{i_1}^{\varepsilon_1} \ldots x_{i_{j-1}}^{\varepsilon_{j-1}}
    + \sum_{1\le j\le t, ~ i_j = k,~ \varepsilon_j=-1} x_{i_1}^{\varepsilon_1} \ldots x_{i_{j}}^{\varepsilon_{j}} = $$
    $$ - \sum_{s+1\le j\le t, ~i_j = k,~ \varepsilon_j=1} x_{i_1}^{\varepsilon_1} \ldots x_{i_{j-1}}^{\varepsilon_{j-1}}
    + \sum_{s+1\le j\le t, ~ i_j = k,~ \varepsilon_j=-1} x_{i_1}^{\varepsilon_1} \ldots x_{i_{j}}^{\varepsilon_{j}} = $$
    \begin{equation}
    \label{eq:D-w}
     - \sum_{s+1\le j\le t, ~i_j = k,~ \varepsilon_j=1} w_{j-1}
    + \sum_{s+1\le j\le t, ~ i_j = k,~ \varepsilon_j=-1}  w_j
    \end{equation}
 Clearly, given $w$ and $s,$ as above one can compute the formal expression (\ref{eq:D-w}) in time  $O(|w|)$. To check if (\ref{eq:D-w}), viewed as an element in $\MZ S_{r,d-1}$, is equal to $0$ it suffices to represent it in the standard group ring form $\sum_{g \in S_{r,d-1}} mg$ (where $m \in \MZ$)  and verify if all coefficients  in this representation are zeros. Now we describe a particular procedure, termed Collecting Similar Terms Algorithm,  which gives the standard group ring form for (\ref{eq:D-w}).
    Given (\ref{eq:D-w}) one can compute in time $O(|w|)$ the following sum
 \begin{equation}
 \label{eq:D-formal}
  - \sum_{s+1\le j\le t, ~i_j = k,~ \varepsilon_j=1} w_{P(j-1)}
    + \sum_{s+1\le j\le t, ~ i_j = k,~ \varepsilon_j=-1} w_{P(j)},
    \end{equation}
Observe now, that two summands $w_p$ and $w_q$ in (\ref{eq:D-formal}) are equal in $S_{r,d-1}$ if and only if $p = q$.  It is easy to see that it takes time $O(|w|)$ to collect similar terms in (\ref{eq:D-formal}), i.e., to compute the coefficients in the standard group ring presentation of (\ref{eq:D-formal}).

It follows that  $\partial w_s/\partial x_k  =  \partial w_t/\partial x_k$ in $\MZ S_{r,d-1}$ if and only if all the coefficients in the standard group ring form of (\ref{eq:D-formal}) are equal to  $0$.   The argument above shows that one can check whether or not $\partial w_s/\partial x_k  =  \partial w_t/\partial x_k$ in $\MZ S_{r,d-1}$ in
 time $O(|w|)$. Since we need to compare all partial derivatives $\partial/\partial x_k, k = 1, \ldots, r,$ of the elements $w_s$ and $w_t$ it takes altogether $O(r|w|)$ time to verify if $w_s = w_t$ in $S_{r,d-1}$.

The routine above allows one to construct effectively the $S_{r,d}$-partition $P_d$ of $D_w$ when given the $S_{r,d-1}$-partition $P_{d-1}$.  A more formal description of the algorithm is given below.

\begin{algorithm}\label{al:partition}{\bf(Computing the $S_{r,d}$-partition of $D_w$)}
    \\{\sc Input.}
Positive integers $r \geq 2$,  $d  \ge 2$, and a word $w \in F(X)$.
    \\{\sc Output.} The $S_{r,d}$-partition function $P_d$  for the set $D_w$.
    \\{\sc Initialization.}
Compute the set $D_w$ and form the initial (trivial) $F$-partition $P_0$ of $D_w$, so $P_0(i) = i$ for every $i \in \{0, \ldots, n\}$.
    \\{\sc Computations.}
\begin{itemize}
    \item[A.]
Compute $A$-partition $P_1$ of $D_w$.
    \item[B.]
For $c = 2 ,\ldots, d$ do:
\begin{itemize}
    \item[1)]
For each $0\le s < t \le n$ such that $w_s = w_t$ in $S_{r,c-1}$
check whether or not $w_s = w_t$ in $S_{r,c}$.
    \item[2)]
Form the $S_{r,c}$-partition $P_d$ of $D_w$.

\end{itemize}
    \item[C.]
Output $P_d$.
\end{itemize}
\end{algorithm}

\begin{lemma}
\label{le:complexity-partition}
Given integers $r, d \geq 2$ and $w \in F$  Algorithm \ref{al:partition}  computes the $S_{r,d}$-partition (the function $P_d$) of $D_w$ in time  $O(d r |w|^3)$.
\end{lemma}

\begin{proof}
Algorithm \ref{al:partition} makes precisely  $d$ iterations $c = 1,
\ldots, d$ by consequently computing the
$S_{r,c}$-partitions  of $D_w$.  After the $S_{r,c-1}$-partition of $D_w$ is computed, the
algorithm computes the  $S_{r,c}$-partition of $D_w$ by comparing elements $w_s, w_t \in D_w$ in $S_{r,c}$. It requires at most $|w|(|w|+1)/2$ such checks, and, as was explained above, every such check can be done in time $O(r|w|)$. Altogether one needs $O(r|w|^3)$ steps to construct the function $P_c$ on the iteration $c$. Since the algorithm makes altogether $d$ iterations it takes it $O(dr|w|^3)$ time to produce $P_d$, as claimed.

\end{proof}

Now we are in a position to show two applications of  Algorithm  \ref{al:partition}.
The first one concerns with computing Fox derivatives in $\MZ S_{r,d}$.

\begin{algorithm}\label{al:derivatives}{\bf(Computing Fox derivatives)}
    \\{\sc Input.}
Positive integers $r \geq 2$,  $d  \ge 2$, a word $w \in F(X)$, and a number $k \in \{1, \dots, r\}$.
    \\{\sc Output.} The standard group ring presentation of the Fox derivative $\partial w/\partial x_k$ in $\MZ S_{r,d}$.
     \\{\sc Computations.}
\begin{itemize}
    \item[A.]
Compute, using formals (\ref{eq:derivative}), the Fox derivative $\partial w/\partial x_k$ in $\MZ F$.
   \item[B.] Compute, using Algorithm  \ref{al:partition}, the $S_{r,d}$-partition of $D_w$.
   \item [C.]  Compute, using the Collecting Similar Terms    Algorithm, the standard group ring form of $\partial w/\partial x_k$ in $\MZ S_{r,d}$

    \item[D.]
Output $\partial w/\partial x_k$ computed in [C.].
\end{itemize}
\end{algorithm}

\begin{lemma}
\label{le:compulexity-derivatives}
Given integers $r, d \geq 2$, a word $w \in F$, and a number $k \in \{1, \dots, r\}$ Algorithm \ref{al:derivatives}  computes the standard group ring presentation of the Fox derivative $\partial w/\partial x_k$ in $\MZ S_{r,d}$  in time  $O(d r |w|^3)$.
\end{lemma}
  \begin{proof}
  Follows from Lemma \ref{le:complexity-partition}.
    \end{proof}
  The second application of   Algorithm  \ref{al:partition} is to WP in $S_{r,d}$.

\begin{algorithm}\label{al:wp_solvable}{\bf(WP in $S_{r,d}$)}
    \\{\sc Input.}
Positive integers $r \geq 2$,  $d  \ge 2$, and a word $w \in F(X)$.
    \\{\sc Output.}
$True$ if $w = 1$ in  $S_{r,d}$ and $False$
otherwise.
    \\{\sc Computations.}
\begin{itemize}
    \item[A.] Compute the set $D_w$.
    \item[B.]
Using Algorithm \ref{al:partition} compute $S_{r,d}$-partition $P_{d}$ of $D_w$.
      \item[C.] If $P_d(0)  = P_d(n)$, i.e., $1 = w$ in $S_{r,d}$,  then output $True$.
      Otherwise output $False$.
\end{itemize}
\end{algorithm}

\begin{theorem}
Algorithm \ref{al:wp_solvable} solves the Word Problem in a free
solvable group $S_{r,d}$ in time $O(d r |w|^3)$.
\end{theorem}

\begin{proof}
Follows immediately from Lemma \ref{le:complexity-partition}.
\end{proof}

\section{Geodesics in  free metabelian groups}
\label{se:metabelian_geodesic}

In this section we discuss the computational hardness of different
variations of Geodesic problems and prove the main result
about NP-completeness of BGLP in free metabelian groups.

\subsection{Algorithmic problems with geodesics in groups}
\label{subsec:geodesics-in-groups}

Let $G$ be a group with a finite   set of generators $X = \{x_1, \ldots,
x_r\}$ and $\mu:F(X) \rightarrow G$ the canonical epimorphism. In this section we view the free group $F(X)$ as the set of all freely reduced words in the alphabet $X^{\pm 1} = X \cup X^{-1}$ with the obvious multiplication.

For a word $w$ in the alphabet $X^{\pm 1}$ by $|w|$ we denote the length of $w$.  The {\em geodesic length} of
an element $g \in G$ relative to $X$, denoted  by $l_X(g)$, is the length of a
shortest word $w \in F(X)$ representing $g$, i.e.,
$$
l_X(g) = \min \{|w| \mid w \in F(X), w^\mu = g\}.
$$
To simplify notation we write, sometimes,  $l_X(w)$ instead of $l_X(w^\mu)$.  A word $w \in F(X)$ is called {\em geodesic} in $G$ relative to $X$, if $|w| = l_X(w)$.

 We are interested here in the following algorithmic {\em search} problem in a given group  $G$ described as above.

\bigskip
\noindent
 {\bf The Geodesic Problem (GP):}

\noindent Given a word $w \in F(X)$ find a geodesic (in $G$) word $\tilde w \in F(X)$ such that $w^\mu = \tilde{w}^\mu$.

\bigskip
One  can consider the following variation of GP.

 \bigskip
\noindent
 {\bf The Geodesic Length Problem (GLP):}

\noindent Given a word $w \in F(X)$ find $l_X(w)$.

\bigskip
Though GLP seems easier than GP (since a solution to GP  gives, in linear time,  a solution to GLP), in practice, to solve GP one usually solves GP first, and only then computes the geodesic length.

As customary in complexity theory  one can modify the search problem GLP to get the corresponding bounded  decision problem:

 \bigskip
\noindent
 {\bf   The Bounded Geodesic Length Problem (BGLP):}

\noindent Let $G$ be a group with a finite generating set $X$. Given a word $w \in F(X)$ and  a natural number $k$ determine if $l_X(w) \leq k$.

\bigskip
It is instructive to compare the algorithmic "hardness" of the problems above and the Word Problem (WP).  Clearly, if one of them is decidable then all of them are decidable. To see the difference we need to recall  a few definitions.  Let $A$ and $B$ be algorithmic problems with inputs sets $I_A$ and $I_B$. Then $A$ is termed {\em Turing reducible} to $B$ in polynomial time, if there exists an algorithm $\CA$  with an oracle for $B$ (which  can be viewed as a "subroutine" of $\CA$  that for a given input $e  \in I_B$ in one step returns the answer for $B$ on $e$) that solves $A$ in polynomial time. Similarly, one can define Turing reducibility in {\em exponential} time. In these cases we write $A \preceq_{T,p} B$ or, correspondingly,  $A \preceq_{T,exp} B$.

Again, it is not hard to see  $WP \preceq_{T,p} BGLP \preceq_{T,p} GLP \preceq_{T,p} GP$. Moreover, since (by brute force algorithm) $GP \preceq_{T,exp} WP$ it follows that all these problems are Turing reducible to each other in exponential time. Moreover, if $G$ has polynomial {\em growth}, i.e., there is a polynomial $p(n)$ such that for each $n \MN$ cardinality of the ball $B_n$  of radius $n$ in the Cayley graph $\Gamma(G,X)$ is at most $p(n)$, then one can easily construct this ball $B_n$ in polynomial time with an oracle for the WP in $G$ (see, for example,  \cite{DLS} for details). It follows that if a group with polynomial growth has WP decidable in polynomial time then all the problems above  have polynomial time complexity with respect to any finite generating set (since the growth and WP stay polynomial for any finite set of generators). Observe now, that by Gromov's theorem \cite{Gromov_pgrowth:1981} finitely generated groups of polynomial growth are virtually nilpotent. It is  also  known that the latter have WP decidable in polynomial time (nilpotent finitely generated groups are linear). These two facts together imply that  the Geodesic Problem is polynomial time decidable in finitely generated virtually nilpotent groups.

On the other hand, there are many groups of exponential growth where GP is decidable in polynomial time, for example, hyperbolic  groups \cite{Epstein}.  Among metabelian groups  the Baumslag-Solitar group $BS(1,2) =  \langle a,t \mid t^{-1} a t = a^2\rangle$  has exponential growth (it is solvable but not polycyclic - and the claim follows from the Milnor theorem \cite{Milnor}) and GP in $BS(1,2)$ is decidable in polynomial time (see \cite{Elder_BS:2007}).

In general, if WP in $G$ is polynomially decidable then $BGLP$ is in the class  NP, i.e., it is decidable in polynomial time by a non-deterministic Turing machine. Indeed, if $l_X(w) \leq k$ then there is  a word $u \in F(X)$ of length at most $k$  which is equal to $w$ in $G$ - this $u$ is a "witness" of polynomial size which allows one to verify  in polynomial time that $l_X(w) \leq k$ (just checking that $u = w$ in $G$). In this case GLP is Turing reducible in polynomial time to an NP problem, but we cannot claim the same for GP. Observe, that BGLP is in NP for any finitely generated metabelian group, since they have WP decidable in polynomial time (see Introduction).

It might happen though,  that WP in a group $G$ is polynomial time decidable, but BGLP in  $G$ is as hard as any problem in the class NP, i.e., it is NP-complete. Recall  (in the notation above) that a decision  problem $B$ is NP-complete if it is in NP and for any decision problem $A$ from NP there is  a  computable in polynomial time function $f:I_A \rightarrow I_B$ (Karp reduction, or a polynomial reduction), such that $A$ is true on $x \in I_A$ if and only if $B$ is true on $f(x)$. The simplest  example of this type is due to   Perry, who  showed in \cite{Parry:1992} that BGLP  is NP-complete in the metabelian group $\MZ_2 wr (\MZ \times \MZ)$ (the wreath product of $\MZ_2$ and $\MZ \times \MZ$).
 In this event, the search  problems GP and GLP are called {\em NP-hard}, this means precisely that some (any) NP-complete problem is a Turing reducible to them in polynomial time.

It would be very interesting to classify finitely generated metabelian groups with respect to computational hardness of their GP or GLP  problems. In the next section we clarify the situation with free metabelian groups. Some remaining open problems are discussed in Section \ref{se:open-problems}.

It was claimed in \cite{DLS} that in free solvable groups of finite rank GLP is decidable in  polynomial time. Unfortunately, in this particular case  their  argument is  fallacious.
Our main result of this section is the following theorem.

\begin{theorem}
\label{th:main} {\bf(Main Theorem)}
Let $M_r$ be a free metabelian group $M_r$ of finite rank $r \geq 2$. Then BGLP in $M_r$ (relative to the standard basis) is NP-complete.
\end{theorem}
\begin{proof}
 The proof of this result consists of two parts. Firstly, in Section \ref{subsec:reduction-to_M2} (Corollary \ref{co:NP-M-2-reduction}) we show that it suffices to prove that BGLP is NP-complete in $M_2$.  Secondly, in Section \ref{subsec:MT} (Theorem \ref{th:NP-complete-M2}) we give a proof that BGLP is, indeed, NP-complete in $M_2$.
\end{proof}

This immediately implies the following results.
\begin{corollary}
The search problems GP and GLP are NP-hard in non-abelian $M_r$ (relative to the standard basis).
\end{corollary}

To prove the Main Theorem we reduce the problem to
the case $r = 2$ and then show that BGLP in $M_2$ is NP-complete.
To see the latter we construct a polynomial reduction
of the Rectilinear Steiner Tree Problem  to BGLP in $M_r$.

\subsection{Reduction to $M_2$}
 \label{subsec:reduction-to_M2}

Let $\CV$  be a variety of groups. For groups $A, B \in \CV$ we denote by $A \ast_\CV B$ the free product of $A$ and $B$ relative to $\CV$. In particular, if $A  = \langle X \mid R\rangle$, and  $B  = \langle Y \mid S\rangle$ are presentations of $A$ and $B$ in $\CV$ then $A \ast_\CV B = \langle x \cup  Y \mid R \cup S\rangle$ is a presentation of $A \ast_\CV B$ in $\CV$.  As usual, $A \ast_\CV B$ satisfies the canonical universal property: any two homomorphism $A \rightarrow C, B \rightarrow C$ into a group $C \in \CV$ extends to a unique homomorphism $A \ast_\CV B \rightarrow C$ (we refer to \cite{Newmann:1967} for details). It follows that if $F_\CV(X \cup Y)$ is a free group in $\CV$ with basis $X \cup Y$ then $F_\CV(X \cup Y) = F_\CV(X) \ast_\CV F_\CV(Y)$.

The following lemma claims  that free $\CV$-factors of a group $G$ are isometrically embedded into $G$.

\begin{lemma}
\label{le:reduction-free-factors}
 Let $A, B \in \CV$ with finite generating sets $X$ and $Y$. Then in the group $A \ast_\CV B$ no geodesic word (relative to $X \cup Y$) for an element from $A$ contains a letter from $Y$. In particular, for any word $w \in F(X)$ its geodesic length in $A$ (relative to $X$) is equal to the geodesic length in $A \ast_\CV B$ (relative to $X  \cup Y$).
\end{lemma}
   \begin{proof}
   Let $w \in F(X)$ be a geodesic word in $A$ relative to $X$. Suppose that  $u \in F(X \cup Y)$ is a geodesic word in $G = A \ast_\CV B$ (relative to $X \cup Y$) that defines the same element as $w$. The identical map $A \rightarrow A$ and the trivial map $B \rightarrow 1$ give rise to a homomorphism  $\phi: G \rightarrow A$. This $\phi$, when applied to $u$, just  "erases" all letters from $Y$. It follows that if $u$ contains a letter from $Y$ then $|u^\phi|< |u| \leq |w|$ - contradiction with the assumption that $w$ is geodesic in $A$ relative to $X$ (since $w= u^\phi$ in $A$).
   \end{proof}

  \begin{corollary}
  \label{co:reduct-to factors2}
   In the notation above, each of the problems GP, GLP,
   BGLP in $A$ (relative to $X$) is polynomial time  reducible to the problem of the same type in $A \ast_\CV B$ (relative to $X \cup Y$).
   \end{corollary}

Notice, that $M_{n+m} = M_n \ast_{\CA_2} M_m$, where $\CA_2$ is the variety of all metabelian groups. Since  WP is  in P for groups from $\CM$ Corollary \label{co:reduct-to factors} implies the following result.
   \begin{corollary}
   \label{co:NP-M-2-reduction}
   If BGLP is NP-complete in $M_2$ then it is NP-complete in  $M_r, r \geq 2,$ relative to the standard bases.
   \end{corollary}

   \begin{remark}
    Corollary \ref{co:NP-M-2-reduction} easily generalizes to free groups in an arbitrary  variety, provided they have WP decidable in polynomial time.
   \end{remark}

\subsection{Rectilinear Steiner Tree Problem}
\label{subsec:RSTP}

The Steiner Tree Problem (STP), which was originally introduced by Gauss, is one of the initial twenty one NP-complete problems that appeared in the Karp's list \cite{Karp:1972}.  We need the following  {\it rectilinear}  variation of STP.

Let $\MR^2$ be the Euclidean plane and $\Gamma$ the integer grid   canonically embedded into $\MR^2$ (all vertices from $\MZ^2$ together with all the horizontal and vertical lines connecting them).  If $A$ is a finite subset of  $\MZ^2$ then a {\em
rectilinear Steiner tree} (RST) for $A$ is a subgraph $T$ of $\Gamma$ such that $A \cup T$ is connected in $\Gamma$. By $s(T)$ (size of $T$) we denote the number of edges in $T$.  An RST for $A$ is  {\em optimal} if it has a smallest possible size among all RST for $A$, we denote such RST by $T_A$. Observe, that a given $A$ may have several  different optimal RST,  but their  size is the same, we denote it by $s(A)$.

Notice that, in general, $T_A$ for $A$ is not a  spanning tree for $A$ in $\Gamma$, it may contain some vertices from $\MZ^2$ which are not in $A$ (so-called, {\em Steiner} points). The Rectilinear Steiner Tree Problem (RSTP) asks for a given finite  $A \subseteq \MZ^2$ and $k \in \MN$ decide if there exists some $T_A$ for $A$ with $s(T_A) < k$.  It is known that RSTP is NP-complete \cite{GJ_rstp1977}.

\subsection{NP-completeness of BGLP in $M_2$}
\label{subsec:MT}

Now we construct a polynomial reduction of RSTP to GLDP in $M_2$ relative to the standard
basis $X = \{x,y\}$.

With each point
$(s,t) \in \MZ^2$ we associate a word
$$w_{s,t} = x_1^s x_2^t \cdot
(x_2 x_1 x_2^{-1} x_1^{-1}) \cdot x_2^{-t} x_1^{-s}$$
in $F(x,y)$. Similarly, with a set of points
 $A = \{(s_1,t_1), \ldots, (s_n,t_n)\} \subset \MZ^2$, ordered in an arbitrary way, we associate a
word
$$w_A = \prod_{i=1}^n w_{s_i,t_i}.$$
Observe, that the word $w_{s,t}$, as well as $w_A$, belongs to $F^\prime = [F,F]$, so in $M_2$ they define elements from $M_2^\prime$.  In particular, the path $p_{w_A}$ is a  closed path in the grid $\Gamma = \MZ^2$, which is viewed as the Cayley graph of the abelianization $F/F^\prime$.

   For  $A \subset \MZ^2$, $(p,q) \in \MZ^2$,  and $m \in \MZ$ we put
   \begin{itemize}
   \item  $A + (p,q) = \{(s+p,t+q) \mid (s,t) \in A\},$
   \item $mA = \{(ms,mt) \mid (s,t) \in A\}.$
   \end{itemize}

\begin{figure}[htbp]
\centering
\epsfig{file=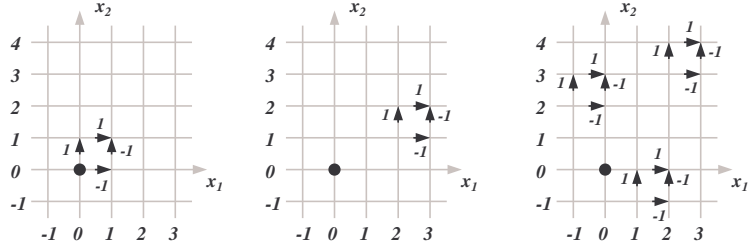, scale=0.6}
\caption{\label{fi:diagram2}
Flows on $\MZ^2$ defined by words $w_{0,0}$, $w_{2,1}$, and
$w_{\{(-1,2),(1,-1),(2,3)\}}$.}
\end{figure}

The following result is due to Hanan \cite{Hanan:1966}.

\begin{theorem} \cite[Theorem 4]{Hanan:1966}
\label{th:Hanan}

Let $A = \{(x_1,y_1), \ldots, (x_n,y_n)\}$ be a finite subset of $\MZ \times \MZ$.
There exists some $T_A$ for $A$ with the set
of Steiner points
    $Q = \{(a_1,b_1), \ldots,(a_q,b_q)\}$
such that $\{a_1, \ldots,a_q\} \subseteq \{x_1,\ldots,x_n\}$ and
$\{b_1, \ldots,b_q\} \subseteq \{y_1,\ldots,y_n\}$.
\end{theorem}

\begin{corollary}\label{co:Hanan}
Let $A$ be a finite subset of $\MZ \times \MZ$. Then
\begin{enumerate}
\item $s(A + (b,c)) = s(A)$ for any $(b,c) \in \MZ^2$;
\item $s(A) = ms(mA)$ for any $m\in \MN$.
\end{enumerate}
\end{corollary}

\begin{proof}
The first statement is obvious because the parallel shift  $(x,y) \rightarrow (x,y) + (b,c)$ is an isomorphism of the Cayley graph $\Gamma$ (in particular, an isometry).

To prove the second statement, notice first that  $s(mA) \le ms(A)$. Indeed, stretching $T_A$ by the factor of $m$ along  all horizontal and vertical lines gives some RST for $mA$, hence the claim.

On the other hand, $s(A) < s(mA)/m$. To see this, observe that by Theorem \ref{th:Hanan}  there exists $R = T_{mA}$ which lies inside the  grid
$m\MZ \times m\MZ$.  Since the coordinates of all vertices in $R$ are
multiples of $k$ one can shrink $R$ by the factor of $m$, in such a way that  the image of $R$ becomes an RST for $A$. Clearly, the size of the image is equal to $s(T_{mA})/m$, hence the result.
\end{proof}

\begin{proposition} \label{pr:geodesic_length}
Let $A$ be a finite subset of $\MZ^2$, $(b,c) \in A$,  and  $n = |A|$. Put $A^\ast = 10n(A-(b,c))$.  Then $l_X(w_{A^\ast}) \in [20ns(A),20ns(A)+4n]$.
\end{proposition}

\begin{proof}
Let $u$ be a geodesic word for $w_{A^\ast}$ relative to the basis $X$. Since $w_{A^\ast} \in F^\prime$ the paths $p_u$ and $p_{w_{A^\ast}}$ are closed paths in $\Gamma =  \MZ^2$ (viewed as the Cayley graph of $M_2/M_2^\prime$). Hence $u$ and $w_{A^\ast}$ determine  the same circulations  $\pi_u =\pi_{A_{w^\ast}}$  on $\Gamma$.  As described in Section \ref{subsec:geodesics},   the flow $\pi_u$ is associated with  the subgraph $\Gamma_u$ generated in $\Gamma$ by $supp(\pi_u) \cup \{(0,0)\}$. It follows from the construction of  the word $w_{A^\ast}$ that the connected components of $\Gamma_u$ are precisely the  $1\times 1$-squares in $\Gamma$, whose lower-left corners are located at the points from $A^\ast$. Notice that $(0,0) \in A -(b,c)$ hence $(0,0) \in A^\ast$.    Now, if $Q$ is a minimal forest for $u$ (a subgraph of $\Gamma$ of minimal size that makes the graph $\Gamma_u \cup Q$ connected in $\Gamma$) then by Theorem \ref{th:geodesics}
\begin{equation}
\label{eq:forest}
|u| = l_X(w_{A^\ast}) = \sum_{e \in supp(p_u)} \pi_u(e) + 2|E(Q)| = 4n + 2|E(Q)|.
\end{equation}
Observe, that an optimal RST $T_{A^\ast}$ for $A^\ast$ also makes the graph $\Gamma_u \cup Q$ connected in $\Gamma$, hence $|E(Q)| \leq s(A^\ast)$. Therefore, $l_X(w_{A^\ast}) \leq 20ns(A)+4n$.

On the other hand assume that $|u| = l_X(w_{A^\ast}) < 20ns(A)$. Hence, from  (\ref{eq:forest}), there exists a minimal forest $Q$ for $A^\ast$ such that $2|E(Q)| < 20ns(A) - 4n$. Since every connected component has precisely $4$ edges and there are $n$ such components, it follows that there is  an RST for $A^\ast$ of size strictly less than
$10ns(A)$ - contradiction with  Corollary \ref{co:Hanan}. This proves the proposition.

\end{proof}

\begin{corollary}
  \label{co:reduction}
  Let $A$ be a finite subset of $\MZ^2$ and $k \in \MN$. In the notation above,
  $$s(A) <  k \Longleftrightarrow l_X(w_{A^\ast}) < 20nk+4n.$$
  In particular, this gives a polynomial reduction of RSTP to BGLP in $M_2$ relative to $X$. \end{corollary}

\begin{proof}
Indeed, if $s(A) < k$ then by Proposition \ref{pr:geodesic_length}  $l_X(w_{A^\ast}) \leq 20ns(A)+4n <  20nk +4n$. On the other hand, suppose $s(A) \geq k$, say $s(A) = k+l$ for some positive $l \in \MN$. Then, again by Proposition \ref{pr:geodesic_length} $l_X(w_{A^\ast}) \geq 20ns(A) =  20n(k+l) > 20nk +4n$, as required.
\end{proof}

\begin{theorem}
\label{th:NP-complete-M2}
GLDP in a free metabelian group $M_2$ is NP-complete.
\end{theorem}

\begin{proof}
 Corollary \ref{co:reduction} gives a polynomial reduction of RSTP in $\MZ^2$ to BGLP in $M_2$. Therefore BGLP in $M_2$ is NP-hard. Meanwhile, as was mentioned above BGLP for $M_2$ is in NP, since WP in $M_2$ is polynomial.
\end{proof}

\section{Open Problems}

\label{se:open-problems}

Denote by $\CM$ the class of all finitely generated metabelian groups.

\begin{problem}
Describe groups in $\CM$ with GP in P.
 In particular, the following partial questions are of interest here:
\begin{itemize}
  \item Are there any groups in $\CM$ with GP  not in P?
  \item  Do polycyclic groups from $\CM$  have GP in P?
   \item Describe wreath products of two f.g. abelian groups have GP in P?
\end{itemize}
\end{problem}

Notice, since WP is in P for groups from $\CM$, it follows that GLP is at most in NP (Turing reducible in P time to BGLP which is in NP). This makes the  following problem very interesting.

\begin{problem}
Classify groups in $\CM$ with NP-complete BGLP.   In particular, the following partial questions are of interest:
\begin{itemize}
  \item  Do polycyclic groups from $\CM$  have GLP in P?
   \item Is it true that a wreath product of finitely generated abelian groups $A wr B$ has NP-complete  BGLP if $A \neq 1$ and the torsion-free rank of $B$ is at least 2?
\end{itemize}
\end{problem}

Clearly,   GLP is polynomial time reducible to GP. On the other hand, it is not clear if there are finitely presented (or finitely generated) groups where GP is not polynomial time reducible to GLP. It would be interesting to clarify the situation in the class $\CM$.  To this end we post the following

\begin{problem}
Are there group is $\CM$ where GP is not polinomial time reducible to GLP?
\end{problem}

\bibliography{../../../main_bibliography}

\end{document}